\documentclass[12pt]{amsart}

\DeclareMathOperator\N{\mathbb N}

\DeclareMathOperator\C{\mathbb C}
\DeclareMathOperator\Z{\mathbb Z}
\DeclareMathOperator\Q{\mathbb Q}
\DeclareMathOperator\Gr{{\mathrm{Gr}}}
\DeclareMathOperator\Fl{\mathbf F}

\DeclareMathOperator{\spa}{span}
\DeclareMathOperator\I{\mathcal I}
\DeclareMathOperator\U{\mathcal U}
\DeclareMathOperator\disc{{\mathrm{D}}}
\DeclareMathOperator\lla{{\boldsymbol{\lambda}}}
\DeclareMathOperator\zz{{\boldsymbol{z}}}
\DeclareMathOperator\CB{CB}
\DeclareMathOperator\SV{SV}
\DeclareMathOperator\CC{\mathcal C}
\DeclareMathOperator\A{\mathcal A}

\newcommand{\glm}{\mathfrak{gl}_m}
\newcommand{\slm}{\mathfrak{sl}_m}

\newtheorem{fact}{Fact}[section]
\newtheorem{lemma}[fact]{Lemma}
\newtheorem{theorem}[fact]{Theorem}
\newtheorem{definition}[fact]{Definition}
\newtheorem{example}[fact]{Example}
\newtheorem{rremark}[fact]{Remark}\newenvironment{remark}{\begin{rremark} \rm}{\end{rremark}}
\newtheorem{proposition}[fact]{Proposition}
\newtheorem{corollary}[fact]{Corollary}
\newtheorem{conjecture}[fact]{Conjecture}

\title[Conformal blocks in the tensor product
of vector representations]
{Conformal blocks in the tensor product
 of vector representations
 and localization formulas}

\author{R. Rim\'anyi}
\address{Department of Mathematics, University of North Carolina at Chapel Hill, USA}
\email{rimanyi@email.unc.edu}

\author{A. Varchenko}
\address{Department of Mathematics, University of North Carolina at Chapel Hill, USA}
\email{anv@email.unc.edu}

\thanks{
The first author is supported by the Marie Curie Fellowship PIEF-GA-2009-235437. He also thanks the hospitality of MPIM Bonn in November 2009. The second author is supported  in part by NSF grant DMS-0555327.}

\begin{document}

\begin{abstract}
Using equivariant localization formulas we give a formula for
conformal blocks  at level one on the sphere as suitable
polynomials. Using this presentation we give a generating set in the
space of conformal blocks at any level if the marked points on the
sphere are generic.

\end{abstract}

\maketitle

\section{Introduction}
We consider conformal blocks on the Riemann sphere
in the $\slm$
Wess-Zimono-Novikov-Witten conformal field theory.  For a partition
$\lla=(\lambda_1,\lambda_2,\ldots,\lambda_m)\in\N^m$ we denote
$|\lla|=\sum \lambda_i$, and $d(\lla)=\lambda_1-\lambda_m$. We fix the
level $\ell$ of the theory with $\ell \geq d(\lla)$ and distinct
points $z_1,\dots,z_{|\lla|},\,\infty$ on the sphere. We
assign to each finite point $z_a$ the standard $m$-dimensional vector
representation of $\slm$, denoted by $V$, and to infinity --
the irreducible $\slm$ representation with highest weight
$(-\lambda_m,\dots, -\lambda_1)$.  The associated space of conformal
blocks
$\CB^\ell_{\zz}(\lla)$
can be realized as a vector
subspace of the tensor
product $V^{\otimes |\lla|}$, and the tensor product  can be realized as a
suitable vector space of polynomials. The subspace of conformal
blocks $\CB^\ell_{\zz}(\lla) \subset V^{\otimes |\lla|}$ is defined as
the set of solutions to a system of differential equations
due to the description of conformal blocks in \cite{FSV1, FSV2}.
We solve that system for $\ell = 1$. In that case
$\dim \CB^1_{\zz}(\lla) = 1$ and we give a formula for one remarkable polynomial
generating the one-dimensional space of conformal blocks
(for $m=2$ a formula was given in \cite{anv_related}).

\medskip

A striking property of the formula
is its similarity to equivariant localization
formulas. According to these formulas, if a torus acts on a compact
manifold with a finite fixed point set, then the integral of an
equivariant cohomology class on the manifold can be computed by
collecting some data at the fixed points. From these data one writes
down a rational function which will be equal to the integral of the
equivariant cohomology class. Not only our conformal block
has the structure of the rational function from a localization
formula, but the proof showing that it is indeed a conformal block
uses equivariant localization. We hope that this connection between
conformal field theory and equivariant cohomology will be useful in
both areas in the future.

\medskip

Taking suitable products of conformal blocks at level one,
we construct
elements in the space of conformal blocks
$\CB^\ell_{\zz}(\lla)$
at any level. We show that
the constructed elements generate $\CB^\ell_{\zz}(\lla)$
for generic $\zz$.
The proofs are based on our formula for conformal blocks at level one.

\medskip

We assign the vector representations to all finite points
$z_1,\dots,z_{|\lla|}$.
The case of more general representations assigned to  finite points
may be studied by fusion procedure.

\medskip

According to a general principle in \cite{MV}, if a space of conformal
blocks is one-dimensional, then there is an associated explicitly calculated
multi-dimensional
Selberg-type integral, giving an integral formula for the
conformal blocks. For $m=2$ such a Selberg-type integral is described
in \cite{anv_related}. We plan to describe the corresponding integral for an
arbitrary $m$ elsewhere.

\medskip

The authors thank V. Schechtman for stimulating and useful discussions.

\medskip

\noindent{\bf Conventions.} The set of natural numbers is
$\N=\{0,1,2,\ldots\}$. Certain expressions in this paper are
parameterized by partitions $\lla$. In notation we put the partition
in bracket, e.g. $\CB(\lla)$, $P(\lla)$, $s(\lla)$. However, for
concrete partitions, e.g. $\lla=(4,2,1)$, we do not repeat brackets,
that is, we do not write $P( \ (4,2,1)\ )$, we will simply write
$P(4,2,1)$.

\section{Spaces of conformal blocks}

\subsection{The Lie algebras $\glm$, $\slm$ and their representations.}
Let $e_{i,j}$, $i,j=1,\dots,m$, be the standard generators of the
complex Lie algebra $\glm$ satisfying the relations
$[e_{i,j},e_{s,k}]=\delta_{j,s}e_{i,k}-\delta_{i,k}e_{s,j}$. We identify the
Lie algebra $\slm$ with the subalgebra in $\glm$ generated by the
elements $e_{i,i}-e_{j,j}$ and $e_{i,j}$ for $i\ne j$, $i,j=1,\dots,m$.

A vector $u$ in a $\glm$-module has weight
$\lla=(\lambda_1,\lambda_2,\dots,\lambda_m)\in\C^m$, if
$e_{i,i}u=\lambda_i u$ for $i=1,\dots,m$.  The vector $u$ is called
{\em singular} if $e_{i,j}u=0$ for $1\leq i<j\leq m$.

Let $\lla=(\lambda_1,\lambda_2,\ldots,\lambda_m)\in\N^m$ be a partition, i.e. let $\lambda_i\geq \lambda_j$ for $i<j$. Throughout the paper we will use the following shorthand notations: $|\lla|=\sum \lambda_i$, and $d(\lla)=\lambda_1-\lambda_m$.
The irreducible finite dimensional $\glm$-module with highest weight $\lla$ will be denoted by $L_{\lla}$. The module $L_{(1,0,\dots,0)}$ is the standard $m$-dimensional vector representation of $\glm$. We will denote it by $V$. The module $V$, and other $\glm$-modules will also be considered as $\slm$-modules.

The Lie algebra $\glm$ acts on $\C[y^{(1)},\dots,y^{(m)}]$ by the differential operators $e_{i,j} \,\mapsto \,y^{(i)} \partial/\partial y^{(j)}$. This action preserves the degree of polynomials. The standard $\glm$-module $V$ is identified with the subspace of $\C[y^{(1)},\dots,y^{(m)}]$ consisting of homogeneous polynomials of degree one.

For a given partition $\lla\in\N^m$, the polynomial ring
\begin{equation}\label{poly_ring}
R_{m,\lla}=\C[y^{(1)}_1,\ldots,y^{(1)}_{|\lla|}, \ \ y^{(2)}_1,\ldots,y^{(2)}_{|\lla|},\ \ \ldots, \ \
y^{(m)}_1,\ldots,y^{(m)}_{|\lla|}]
\end{equation}
admits a $\glm$-module structure by the rule $e_{i,j}\mapsto \sum_{a=1}^{|\lla|} y^{(i)}_a \partial/\partial y^{(j)}_a$.
The $\glm$-module $V^{\otimes |\lla|}$ is identified with the subspace of $R_{m,\lla}$  consisting of polynomials that have homogeneous degree 1 with respect to each $m$-tuple of variables $y^{(1)}_a,\ldots,y^{(m)}_a$, for $a=1,\ldots,|\lla|$.

\subsection{Conformal blocks.}

Consider the Lie algebra $\glm$, its vector representation $V$, a partition $\lla$, and the $\glm$-module $V^{\otimes|\lla|}$ as identified with a subspace of $R_{m,\lla}$ in (\ref{poly_ring}). Recall the action of $e_{i,j}$ on the latter, namely
\begin{equation}\label{eq:eij}
e_{i,j}=\sum_{a=1}^{|\lla|} y^{(i)}_a \partial/\partial y^{(j)}_a.
\end{equation}
The space of singular vectors of weight $\lla$ is
$$\SV(\lla) =  \{p\in V^{\otimes|\lla|}\ |\ e_{i,j}p=0, e_{i,i}p=\lambda_ip \ \text{for}\ 1\leq i<j \leq m\}.$$


\noindent Fixing distinct complex numbers $\zz=(z_1,\ldots,z_{|\lla|})$
we define the differential operators
\begin{equation}\label{eq:ez}
e^{\zz}_{i,j}=\sum_{a=1}^{|\lla|} z_a y_a^{(i)} \frac{\partial}{\partial y^{(j)}_a}, \qquad\text{for}\ 1\leq i,j\leq m.
\end{equation}

\noindent For a positive integer $\ell\geq d(\lla)$ we define the {\em space of conformal blocks at level $\ell$} by
\[
\CB_{\zz}^\ell(\lla)=\{p\in \SV(\lla)\ |\ \left(   e^{\zz}_{1,m}    \right)^{\ell-d(\lla)+1} p=0\}.
\]

If $\ell-d(\lla)+1$ is greater than the $y^{(m)}$-degree (namely,
$\lambda_m$) of singular vectors from $\SV(\lla)$, then the defining
equation of level $\ell$ conformal blocks is vacuous, hence we have
$$
\CB_{\zz}^{d(\lla)}(\lla) \subset \CB_{\zz}^{d(\lla)+1}(\lla)
\subset \ldots\subset \CB_{\zz}^{\lambda_1}(\lla)=\SV(\lla).
$$

Let us emphasize that in the definition above, as well as in the whole
paper, $\zz$ denotes a collection of {\em distinct} complex numbers.

\begin{remark}\label{rem:cft}
This definition of conformal blocks is nonstandard. Usually the space
of conformal blocks in a WZW model is defined if one has a set of distinct points on
a Riemann surface marked with irreducible representations of an affine
Lie algebra, see \cite{KL}. If the Riemann surface is the Riemann
sphere, then one can describe the space of conformal blocks in terms
of finite dimensional representations of the corresponding finite
dimensional Lie algebra. That description is one of two main results
of \cite{FSV1} and \cite{FSV2}.  We take that description as our
definition for the case when the marked points of the Riemann sphere
are $z_1,\dots, z_{|\lla|},\infty$ and the associated representations
are the standard $\slm$-modules assigned to all finite points $z_a$
and the irreducible $\slm$ highest weight module assigned to the point
at infinity with the highest weight  $(-\lambda_m,\dots,-\lambda_1)$.
\end{remark}

\subsection{``Symmetry and vanishing" description of conformal blocks for
$\lla = (N,\dots,N)$}
Assume that $\lla=(N,\dots,N)\in \N^m$. In this case  we can interpret the space of conformal blocks as follows.

Let $\CC=\C^{m|\lla|}$ be the vector space with coordinates
$y^{(i)}_{a}$ for $i=1,\ldots,m$, $a=1,\ldots,|\lla|$. The polynomials
$p\in V^{\otimes |\lla|}$ are functions on $\CC$.  Consider the vector
subspace
\begin{eqnarray*}
\A(\zz) =  \{ \gamma \in \CC\ | \ y^{(m)}_{a}({\gamma})  =  z_a y^{(1)}_{a}(\gamma),\ a=1,\dots,|\lla| \}
\end{eqnarray*}
of codimension $|\lla|$ in $\CC$.

For any integer $1\leq k\leq |\lla|$ and a subset $B_k = \{1\leq b_1<
\dots < b_k\leq |\lla|\}$ introduce the differential operator
\begin{eqnarray*}
\partial_{B_k}\ = \ \prod_{i=1}^k \, y^{(1)}_{b_i}   \frac{\partial}{\partial y^{(m)}_{b_i}}.
\end{eqnarray*}

The special linear {\em group} $SL_m$ acts diagonally on the polynomial ring
$R_{m,\lla}$ of (\ref{poly_ring}) by substitution in each set
$y^{(1)}_a, \ldots, y^{(m)}_a$ of $m$-tuples of variables. Hence it
acts on the subspace identified with $V^{\otimes |\lla|}$ too.
If $\lla = (N,\dots,N)$, then the subspace of singular vectors $SV(\lla)\subset V^{\otimes |\lla|}$ is the subspace of
$SL_m$-invariant polynomials.

The following theorem is the $SL_m$ analogue of \cite[Thm. 4.3]{R}, 
also \cite[Lemma 1.3]{LV}.  Its proof is straightforward calculation.

\begin{theorem} \label{thm vanish of CoB}
Let $\lla=(N,\dots,N)$ for some $N\in \N$. Then an
 $SL_m$-invariant polynomial $p\in V^{\otimes |\lla|}$ lies in $\CB^\ell_{\zz}(\lla)$, if and only if
\begin{eqnarray} \label{vanish formula}
(\partial_{B_k}\,p) |_{\A(\zz)}  =  0
\end{eqnarray}
for all $B_k = \{1\leq b_1< \dots < b_k\leq |\lla|\}$ with $k \leq
N-\ell-1$. In other words, an $SL_m$-invariant polynomial
$p\in V^{\otimes |\lla|}$ lies in $\CB^\ell_{\zz}(\lla)$, the space of
conformal blocks, if and only if it vanishes at $\A(\zz)$ to order
$\geq N-\ell$.
\end{theorem}

\begin{corollary} \label{cor vanish}
If $\lla=(N,\dots,N)$ for some $N\in \N$, then an
 $SL_m$-invariant polynomial $p\in V^{\otimes |\lla|}$ lies in
$\CB^\ell_{\zz}(\lla)$ if and only if $p$ vanishes to order $\geq
N-\ell$ at the $SL_m$-orbit of $\A(\zz)$.
\end{corollary}

\section{Rational functions: localizations, divided differences}\label{section:rational}

\subsection{Rational function identities obtained from localization formulas.}

In this section we collect some identities for rational functions that
will be useful in producing elements in the spaces of conformal blocks.  The
resultant of the sets $A$ and $B$ of variables is defined to be
$$R(A|B)=\prod_{a\in A} \prod_{b\in B} (a-b).$$ Later we will also
drop the set signs, and write eg. $R(a|b,c)$ for
$R(\{a\}|\{b,c\})=(a-b)(a-c)$.  For sets $A_1,\ A_2,\ldots, A_m$ of
variables $R(A_1|A_2|\ldots|A_m)$ we will denote the generalized
resultant $\prod_{i<j} R(A_i|A_j)$. For a set $I$ of indices, $z_I$
will denote the set of variables $\{z_i\}_{i\in I}$.

\begin{lemma}\label{lagrange}
Let $k<n$ be positive integers. Let $p$ (resp. $q$) be symmetric
polynomials in $k$ (resp. $n-k$) variables.  Let
$\binom{\{1,\ldots,n\}}{k}$ denote the set of $k$-element subsets of
$\{1,\ldots,n\}$. For such a subset $I$, let $\bar{I}$ denote the
complement set $\{1,\ldots,n\}-I$. Then for $\deg(p)+\deg(q)<k(n-k)$,
$$\sum_{I\in \binom{\{1,\ldots,n\}}{k}} \frac{p(z_I)q(z_{\bar{I}})}{R(z_I|z_{\bar{I}})}=0$$
is an identity of rational functions in the variables $z_1,\ldots,z_n$.
\end{lemma}

This lemma is well known in at least two areas of mathematics: the
theory of symmetric functions (in terms of generalized Lagrange
interpolation, see e.g. \cite[Thm. 7.7.1]{lascoux}), and in the theory
of equivariant localization. Let us recall this latter argument. If
the torus $T$ acts on the compact manifold $M$ with fixed points
$f_1,\ldots, f_r$, and $\alpha \in H_T^*(M;\Q)$ is an equivariant
cohomology class, then the Atiyah-Bott equivariant localization
theorem \cite{ab} states that
$$\int_M \alpha = \sum_{i=1}^r \frac{\alpha|_{f_i}}{e(T_{f_i}M)}.$$
Here $e(T_fM)$ is the $T$-equivariant Euler class of the
representation of $T$ on the tangent space of $M$ at $f$. Let us
choose $M$ to be the Grassmannian $\Gr_k(\C^n)$, with the $T=(S^1)^n$
action induced by the standard $T$-action on $\C^n$. Let $\alpha$ be
the $p$-value of the Chern roots of the universal subbundle over $M$,
times, the $q$-value of the Chern roots of the universal quotient
bundle over $M$. Since $\deg(p)+\deg(q)<k(n-k)=\dim(\Gr_k\C^n)$ the
integral $\int_M \alpha$ is clearly 0. Then the localization theorem
proves Lemma \ref{lagrange}.

\smallskip

Consider the partition $\lla=(\lambda_1,\ldots,\lambda_m) \in
\N^m$. Correspondingly, we can consider the partial flag manifold
$\Fl_{\lla}$ parameterizing the flags of linear subspaces
$$0=V_0 \subset V_1 \subset V_2 \subset \ldots \subset V_{m-1} \subset
V_m = \C^{|\lla|}$$ in $\C^{|\lla|}$, where $V_j/V_{j-1}$ has
dimension $\lambda_j$ for $j=1,2,\ldots,m$. Let the tautological
bundle over $\Fl_{\lla}$ corresponding to the $j$'th linear space
$V_j$ be called $E_j$ ($j=0,1,\ldots,m$). For $j=1,\ldots,m$ let
$\gamma_j$ be the collection of the Chern roots of the quotient bundle
$E_j/E_{j-1}$. The standard action of $(S^1)^{|\lla|}$ on
$\C^{|\lla|}$ induces an action on $\Fl_{\lla}$. Then for the
symmetric polynomials $p_j$ in $\lambda_j$ variables ($j=1,\ldots,m$),
the equivariant localization theorem gives
\begin{equation}\label{flag}
\int_{\Fl_{\lla}}\ \  \prod_{i=1}^m p_j(\gamma_j) = \sum_{\I} \frac{  \prod_{i=j}^m p_j(z_{I_j})  }
{R(z_{I_{1}}|z_{I_2}|\ldots|z_{I_m})},
\end{equation}
where $\I=(I_1,I_2,\ldots,I_m)$ is a partitioning of the integers
$\{1,\ldots,|\lla|\}$ into $m$ parts satisfying
$$\cup_j I_j=\{1,\ldots,|\lla|\}, \qquad I_i\cap I_j =\emptyset \ (i\not= j),\qquad |I_j|=\lambda_j.$$

When $\sum \deg(p_i)<\dim \Fl_{\lla}=\prod_{i<j}(\lambda_i-\lambda_j)$, then the left hand side of (\ref{flag}) vanishes, establishing the identity that the right hand side is 0.

\subsection{Localization formulas vs divided differences}\label{section:dd}

For a polynomial $p$ in variables $z_1, z_2,\ldots$, the $i$th divided difference is defined by
\begin{equation}\label{dd_def}
\partial_i p = \frac{ p - p(z_i \leftrightarrow z_{i+1}) }{z_i-z_{i+1}}.
\end{equation}
It is well known that the divided difference operators $\partial_i$ satisfy the relations
\begin{equation}\label{eqn:coxeter}
\partial_i\partial_j=\partial_j\partial_i\ \text{if}\ |i-j|\geq2,\qquad\qquad\partial_i \partial_{i+1} \partial_i= \partial_{i+1} \partial_{i} \partial_{i+1}.
\end{equation}
Hence $\partial_{\omega}$ can be defined for a permutation of indexes,
as $\partial_{i_r}\cdots\partial_{i_1}$, where $s_{i_1}\cdots s_{i_r}$ is a reduced word for $\omega$ in terms of the elementary transpositions $s_i$.

The algebra of divided differences is closely related with localization formulas. For example, let $k<n$, and let $\omega_{(k,n-k)}$ be the
following permutation of $\{1,\ldots,n\}$:
$$i\mapsto (n-k)+i \ \text{for}\ i\leq k \qquad\text{and}\qquad  i\mapsto i-k \ \text{for}\  i>k.$$
Let $p$ and $q$ be symmetric polynomials in $k$ and $n-k$ variables respectively. Then
$$\partial_{\omega_{(k,n-k)}}\left( p(z_1,\ldots,z_k)q(z_{k+1},\ldots,z_n)\right) = \sum_{I\in \binom{\{1,\ldots,n\}}{k}} \frac{p(z_I)q(z_{\bar{I}})}{R(z_I|z_{\bar{I}})}.$$

\smallskip

More generally, consider the partition $\lla$ again. We define the permutation $\omega_{\lla}$ as follows. For
$\lambda_1+\ldots+\lambda_{u-1}< i \leq \lambda_1+\ldots+\lambda_{u}$,
$$i \mapsto i+\sum_{j>u} \lambda_j - \sum_{j<u} \lambda_j.$$
[In plain language $\omega_{\lla}$ is described by dividing the numbers from 1 to $|\lla|$ into groups of cardinality $\lambda_i$ in order, then
reversing the order of the groups without changing the relative positions of pairs of numbers in the same group.]

Then for the symmetric polynomials $p_j$ in $\lambda_j$ variables ($j=1,\ldots,m$), we have
$$\partial_{\omega_{\lla}} \left(\prod_{j=1}^m p_j(z_{J_j})\right) = \sum_{\I} \frac{  \prod_{i=j}^m p_j(z_{I_j})  }
{ R(z_{I_1}|z_{I_2}|\ldots|z_{I_m})},$$
where the summation runs for partitions $\I$ as in (\ref{flag}).

\subsection{A generalized divided difference.}\label{section:general_dd}

In the next chapter we will mention an extended version of divided difference operations. This applies to functions depending not only on
$z_1,z_2,\ldots$, but on other sets of variables. We modify only the numerator in the definition (\ref{dd_def}) by
applying the transposition $i \leftrightarrow (i+1)$ to certain sets of variables (including the $z$'s). These sets of variables
will be indicated in the upper index of the $\partial$ sign, eg.
$$\partial^{x,y,z}_i p(x_1,x_2,\ldots, y_1,y_2,\ldots, z_1,z_2,\ldots)=\frac{p-p(z_i \leftrightarrow z_{i+1}, x_i\leftrightarrow x_{i+1}, y_i \leftrightarrow y_{i+1})}{z_i-z_{i+1}}.$$
Generalized divided differences with fixed upper indexes also satisfy the relations (\ref{eqn:coxeter}), hence $\partial^{x,y,z}_\omega$ is defined for permutations $\omega$. Observe, however, that our generalized divided difference operators do not preserve polynomials.

\section{The $P_{\zz}(\lla)$ function}
\label{section:Pmn_properties}

\subsection{Definition, examples}\label{section:def_P}

Let $m\geq 2$, and let $\lla=(\lambda_1,\ldots,\lambda_m)\in \N^m$ be a partition. We will study various expressions in the variables
$$y^{(j)}_a\ \text{and}\ z_a,\quad\text{for}\quad j\in \{1,\ldots,m\}, a\in \{1,\ldots,|\lla|\}.$$
For a subset $U\subset \{1,\ldots,|\lla|\}$, we define $Y^{(j)}_U=\prod_{a\in U}y^{(j)}_a$.

\begin{definition}
We define
$$P_{\zz}(\lla)=\sum_{\I} \frac{ \prod_{j=1}^m Y_{I_j}^{(j)}}{R(z_{I_1}|z_{I_2}|\ldots|z_{I_m})},$$
where the summation runs for $\I=(I_1,\ldots,I_m)$ with $I_i\cap I_j=\emptyset$, $\cup I_j=\{1,\ldots,|\lla|\}$, $|I_j|=\lambda_j$.
\end{definition}

Recall that associated with the partition $\lla$ we defined the
permutation $\omega_{\lla}$ in Section \ref{section:dd}, and recall
the extended divided difference operation from Section
\ref{section:general_dd}. The function $P_{\zz}(\lla)$ can be written
in the concise form
\begin{equation}\label{eq:concise}
P_{\zz}(\lla)=\partial^{y^{(1)},\ldots,y^{(m)},z}_{\omega_{\lla}}
\left( 
y^{(1)}_1 \cdots  y^{(1)}_{\lambda_1}
y^{(2)}_{\lambda_1+1} \cdots y^{(2)}_{\lambda_1+\lambda_2}
\ \ \cdots\ \
y^{(m)}_{|\lla|-\lambda_m+1} \cdots  y^{(m)}_{|\lla|}
\right) .
\end{equation}

\begin{example}\rm
\begin{align*}
P_{\zz}(1,1)=& \frac{ y^{(1)}_1 y^{(2)}_2}{z_1-z_2} + \frac{ y^{(1)}_2 y^{(2)}_1}{z_2-z_1}
= \frac{\det  \begin{pmatrix} y^{(1)}_1 & y^{(1)}_2 \\
                             y^{(2)}_1 & y^{(2)}_2 \end{pmatrix} }{z_1-z_2}.\\
P_{\zz}(2,1)=&\frac{y^{(1)}_1 y^{(1)}_2 y^{(2)}_3 }{(z_1-z_3)(z_2-z_3)} + \frac{y^{(1)}_1 y^{(1)}_3 y^{(2)}_2 }{(z_1-z_2)(z_3-z_2)} +\frac{y^{(1)}_2 y^{(1)}_3 y^{(2)}_1 }{(z_2-z_1)(z_3-z_1)}.\\
P_{\zz}(2,2)=&  \frac{ y^{(1)}_1 y^{(1)}_2 y^{(2)}_3 y^{(2)}_4 }{ R(z_1,z_2|z_3,z_4) } + \text{[5 similar terms]}  \\
=& \frac{\text{perm}\begin{pmatrix}
y^{(1)}_1 y^{(1)}_2 & y^{(1)}_3 y^{(1)}_4 \\
y^{(2)}_1 y^{(2)}_2 & y^{(2)}_3 y^{(2)}_4
\end{pmatrix}}{ R(z_1,z_2|z_3,z_4) } + \frac{\text{perm}\begin{pmatrix}
y^{(1)}_1 y^{(1)}_3 & y^{(1)}_2 y^{(1)}_4 \\
y^{(2)}_1 y^{(2)}_3 & y^{(2)}_2 y^{(2)}_4
\end{pmatrix}}{ R(z_1,z_3|z_2,z_4) } +\frac{\text{perm}\begin{pmatrix}
y^{(1)}_1 y^{(1)}_4 & y^{(1)}_2 y^{(1)}_3 \\
y^{(2)}_1 y^{(2)}_4 & y^{(2)}_2 y^{(2)}_3
\end{pmatrix}}{ R(z_1,z_4|z_2,z_3) },\\
& \text{[where $perm$ means the permanent of a matrix.]} \\
P_{\zz}(1,1,1)=&\sum_{(i,j,k)\in S_3} \frac{ y^{(1)}_i y^{(2)}_j y^{(3)}_k}{(z_i-z_j)(z_i-z_k)(z_j-z_k)}=
\frac{\det\begin{pmatrix} y^{(1)}_1 & y^{(1)}_2 & y^{(1)}_3 \\
y^{(2)}_1 & y^{(2)}_2 & y^{(2)}_3 \\
y^{(3)}_1 & y^{(3)}_2 & y^{(3)}_3
\end{pmatrix}}{(z_1-z_2)(z_1-z_3)(z_2-z_3)}.
\end{align*}
\end{example}

\begin{remark} \label{rmk:sasha_m2}
In \cite[(3.2)]{anv_related} the rational function
$$
\det \left( \frac{y^{(1)}_{j+N} y^{(2)}_{i} - y^{(1)}_i y^{(2)}_{j+N} }{z_i-z_{j+N}} \right)_{i,j=1,\ldots,N}
$$
is considered (with the substitution $y^{(1)}_{j}=1$ for all
$j$). One can show that this function is equal to
$$\pm \left(\prod_{1\leq i<j \leq N} (z_i-z_j)(z_{i+N}-z_{j+N}) \right)\cdot P_{\zz}(N,N).$$
\end{remark}

\subsection{Vanishing properties} \label{section:vanishing}
The function $P_{\zz}(\lla)$ satisfies remarkable differential
equations with respect to variables $y^{(i)}_j$. 
Recall the differential operators 
from (\ref{eq:eij}).


\begin{theorem}\label{thm:symmetry} We have
$$
e_{k,l} P_{\zz}(\lla)=0 \qquad \text{if}\qquad
\lambda_k \geq \lambda_l.$$
\end{theorem}

\begin{proof} The effect of the terms of the differential operator $e_{k,l}$ on a square-free
$y$-monomial is that they
replace a $y^{(l)}$ variable with the corresponding $y^{(k)}$
variable. Hence we have
$$e_{k,l}\left( \prod_{j=1}^m Y_{I_j}^{(j)}\right)
= \sum_{v\in I_l} \left( \prod_{j=1}^m Y_{I_j}^{(j)} \frac{y^{(k)}_v}{y^{(l)}_v}\right).$$
Therefore, the $y$-monomials occurring in $e_{k,l}P_{\zz}(\lla)$ are of the form $\prod_{j=1}^m Y^{(j)}_{K_j}$
for $K_i\cap K_j=\emptyset$, $\cup K_j=\{1,\ldots,|\lla|\}$, and
$$|K_j|=\begin{cases} \lambda_j & j\not=k,l \\
\lambda_k+1 & j=k \\
\lambda_l-1 & j=l.\end{cases}$$
For such a monomial, and a $v\in K_k$ define
$$I^{(v)}_j=\begin{cases} K_j & j\not= k,l \\
K_k-\{v\} & j=k \\
K_l \cup \{v\} & j=l.\end{cases}$$
Then the coefficient of $\prod_{j=1}^m Y^{(j)}_{K_j}$ in $e_{k,l}P_{\zz}(\lla)$ is
\begin{equation}\label{eab_calculation}
\sum_{v\in K_k} \frac{1}{R(z_{I^{(v)}_1}|\ldots|z_{I^{(v)}_m})}= \frac{\pm 1}{R(z_{K_1}|z_{K_2}|\ldots|z_{K_m})}
\sum_{v\in K_k} \frac{R(z_v|z_{K_l})}{R(z_v|z_{K_k-\{v\}})}.
\end{equation}
The main observation in the last equality is that the sign $\pm 1$
does not depend on the choice of $v\in K_k$.

In the factor
$$\sum_{v\in K_k} \frac{R(z_v|z_{K_l})}{R(z_v|z_{K_k-\{v\}})}$$ of the
last expression we consider the variables $z_{K_l}$ as parameters. In
the remaining variables this expression is of the form of the identity
in Lemma \ref{lagrange}. Hence it is 0 if the numerator has smaller
degree than the denominator, ie. if
$\lambda_l-1=|K_l|<|K_k|-1=\lambda_k$. This inequality is satisfied if
$\lambda_k\geq \lambda_l$.
\end{proof}

Now recall the differential operators from (\ref{eq:ez}).


\begin{theorem} \label{thm:e_z}
We have
\begin{align}
\label{ez_1} e^{\zz}_{k,l} P_{\zz}(\lla) =0 & \qquad\text{      if    } \lambda_k>\lambda_l ,\\
\label{ez_2} \left(e^{\zz}_{k,l}\right)^2 P_{\zz}(\lla)=0 & \qquad\text{      if     }  \lambda_k = \lambda_l .
\end{align}
\end{theorem}

\begin{proof} The proof of (\ref{ez_1})
is analogous with the proof of Theorem \ref{thm:symmetry}. The change is that formula~(\ref{eab_calculation})
expressing the coefficients of $y$-monomials of
$e^{\zz}_{k,l}P_{\zz}(\lla)$ has to be replaced with
\begin{equation}\label{eab_calculation1}
\sum_{v\in K_k} \frac{1}{R(z_{I^{(v)}_1}|\ldots|z_{I^{(v)}_m})}= \frac{\pm 1}{R(z_{K_1}|z_{K_2}|\ldots|z_{K_m})}
\sum_{v\in K_k} \frac{z_vR(z_v|z_{K_l})}{R(z_v|z_{K_k-\{v\}})}.
\end{equation}

The last factor $\sum_{v\in K_k}
\frac{z_vR(z_v|z_{K_l})}{R(z_v|z_{K_k-\{v\}})}$ vanishes due to Lemma
\ref{lagrange} if the numerator has smaller degree than the
denominator, i.e. if $1+(\lambda_l-1)=1+|K_l|<|K_k|-1=\lambda_k$. This
holds if $\lambda_k>\lambda_l$.

\smallskip

To prove (\ref{ez_2}) observe that the effect of the terms of the
differential operator $\left(e^{\zz}_{k,l}\right)^2$ on square-free
$y$-monomials is replacing two $y^{(l)}$ variables with the
corresponding $y^{(k)}$ variables times the corresponding $z$
variables.  That is, we have
$$\left(e^{\zz}_{k,l}\right)^2\left( \prod_{j=1}^m Y_{I_j}^{(j)}\right) = \sum_{v\not= w\in I_l} \left( \prod_{j=1}^m Y_{I_j}^{(j)} \frac{y^{(k)}_v y^{(l)}_w z_vz_w}{y^{(k)}_v y^{(l)}_w   }\right).$$
Therefore, the $y$-monomials occurring in $(e^{\zz}_{k,l})^2P_{\zz}(\lla)$ are of the form $\prod_{j=1}^m Y^{(j)}_{K_j}$
for $K_i\cap K_j=\emptyset$, $\cup K_j=\{1,\ldots,|\lla|\}$, and
$$|K_j|=\begin{cases} \lambda_j & j\not=k,l \\
\lambda_k+2 & j=k \\
\lambda_l-2 & j=l.\end{cases}$$
For such a monomial, and $v\not= w\in K_k$ define
$$I^{(v,w)}_j=\begin{cases} K_j & j\not= k,l \\
K_k-\{v,w\} & j=k \\
K_l \cup \{v,w\} & j=l.\end{cases}$$
Then the coefficient of $\prod_{j=1}^m Y^{(j)}_{K_j}$ in $(e^{\zz}_{k,l})^2P_{\zz}(\lla)$ is
$$\sum_{v\not=w\in K_k} \frac{z_vz_w}{R(z_{I^{(v,w)}_1}|\ldots|z_{I^{(v,w)}_m})}= \frac{\pm 1}{R(z_{K_1}|\ldots|z_{K_m})}
\sum_{v\not= w\in K_k} \frac{z_vz_wR(z_v,z_w|z_{K_l})}{R(z_v,z_w|z_{K_k-\{v,w\}})}.$$
The main observation in the last equality is, again, that the sign $\pm 1$ does not depend on the choice of $v\not= w\in K_k$.

The factor
$$\sum_{v\not= w\in K_k} \frac{z_vz_wR(z_v,z_w|z_{K_l})}{R(z_v,z_w|z_{K_k-\{v,w\}})}$$
vanishes---because of Lemma \ref{lagrange}---if the degree $2+2(\lambda_l-2)$ of the numerator is less than the degree $2\lambda_k$ of the denominator, i.e. if $\lambda_k\geq \lambda_l$.
\end{proof}

\subsection{$P_{\zz}(\lla)$ functions are conformal blocks.}

Let $m\in \N$, let $\lla\in \N^m$ be a partition, and let $\ell\geq
d(\lla)$, as before. Recall that a function in the variables
$y^{(j)}_a$, $z_a$ ($j=1,\ldots,m$, $a=1,\ldots, |\lla|$) belongs to
the space $\CB^{\ell}_{\zz}(\lla)$ of conformal block of level $\ell$,
if
\begin{itemize}
\item it is a polynomial of homogeneous degree 1 with respect to each
  $m$-tuple of variables $y^{(1)}_a, \ldots, y^{(m)}_a$,
\item it vanishes under the action of the differential operator
  $\sum_{a=1}^{|\lla|} y^{(k)}_a \partial/\partial y^{(l)}_a$ for all
  $k<l$,
\item it vanishes under the action of the differential operator
  $(\sum_{a=1}^{|\lla|} z_a y^{(1)}_a \partial/\partial
  y^{(m)}_a)^{l-d(\lla)+1}$.
\end{itemize}

The properties of the $P_{\zz}(\lla)$ function presented in Sections
\ref{section:def_P}, \ref{section:vanishing} translate to the
following theorem.

\begin{theorem} \label{thm:p_in_CB}
We have
\begin{align*}
P_{\zz}(\lla)\in \CB_{\zz}^{d(\lla)}(\lla) & & \text{if    } d(\lambda)>0, \\
P_{\zz}(\lla)\in \CB_{\zz}^{1}(\lla)\ \ \   & & \text{if    }  d(\lambda)=0.
\end{align*}
\end{theorem}

Dimensions of $\CB_{\zz}^{\ell}(\lla)$ spaces are independent of
$\zz$. They are calculated, e.g. as structure constants in so-called
fusion rings, see e.g. \cite{Z}, \cite{GV}. A special case of these
results is that $\dim \CB^1_{\zz}(N+1,\ldots,N+1,N,\ldots,N)=1$ for any
$m$ and any number of $N$'s. Hence the relevant part of Theorem \ref{thm:p_in_CB}
 can be reformulated as follows.

 \begin{theorem}  \label{thm:P_basis}
 Let $m\geq 2$ and let
 $\lla=(N+1,\ldots,N+1,N,\ldots,N)\in \N^m$ for any number of $N$'s
 (possibly $m$).
 Then
 for any collection of pairwise distinct complex numbers
 $\zz=(z_1,\ldots,z_{|\lla|})$, the polynomial
 $P_{\zz}(\lla)$ is a basis in the one-dimensional space of conformal
 blocks $\CB_{\zz}^1(\lla)$. \qed
 \end{theorem}

 For $m=2$ and $\lla=(N,N)$ a basis conformal block is given in \cite{anv_related}, cf. Remark \ref{rmk:sasha_m2}.

\subsection{Asymptotics.}
When concrete formulas are known for conformal blocks, one can often
derive Selberg integral formula type applications, see
e.g. \cite{MV}, \cite{FSV}, \cite{anv_related}.
In these applications one needs certain
asymptotic properties of conformal blocks. Now we show such an asymptotic property of
$P_{\zz}(\lla)$.

We define the discriminant of the ordered list of variables $x_1,\ldots,x_n$ by
$$\disc(x_1,\ldots,x_n)=\prod_{1\leq i < j\leq n} (x_i-x_j).$$
Let
$$\check{P}(N)=P_{\zz}(\underbrace{N,N,\ldots,N}_m)\cdot \prod_{k=1}^N \disc(z_{(k-1)m+1},\ldots,z_{km}).$$
Observe that
$$\check{P}(1)=\det\left(y^{(i)}_j\right)_{\genfrac{}{}{0pt}{}{i=1,\ldots,m}{j=1,\ldots,m}}.$$
Direct calculation gives the following theorem.

\begin{theorem}\label{thm:asymptotics}
For $N\geq 2$ the limit of $\check{P}(N)$ as
$$z_{(k-1)m+1}, z_{(k-1)m+2},\ldots,z_{km-1} \to z_{km} \qquad \forall k=1,\ldots,N$$
is
$$-\frac{\prod_{k=1}^N \det\left(y^{(i)}_j\right)
_{\genfrac{}{}{0pt}{}{i=1,\ldots,m\ \ \ \ \ \ \ \ }{j=(k-1)m+1,\ldots,km}}}{\disc(z_m,z_{2m},\ldots,z_{Nm})^{m(m-1)}}.$$
\end{theorem} \qed

\section{$P_{\zz}(N,\dots,N)$ satisfies the KZ equations at level 1.}\label{kz}

Let
$$
X_{|\lla|}\ = \ \{ \zz=(z_1,\dots,z_{|\lla|}) \in
\C^{|\lla|}\ | \ z_a\neq z_b\ {\rm for\ all }\ a\neq b\}\ .
$$
The trivial vector  bundle
$\eta  : V^{\otimes |\lla|}\times X_{\lla}\to X_{\lla}$
has the KZ connection at level $\ell$,
$$
\frac{\partial}{\partial z_i} \ - \ \frac 1{m+\ell}\,
\sum_{j\neq i} \frac {\Omega^{(i,j)}}{z_i-z_j}\ ,
\qquad
i=1,\dots,|\lla|\ ,
$$
where
$\Omega = \oplus_{a,b=1}^m e_{a,b}\otimes e_{b,a}$ is the $\glm$ Casimir operator
and $\Omega^{(i,j)}:V^{\otimes |\lla|}\to V^{\otimes|\lla|}$
is the linear operator acting as $\Omega $ on the
$i$-th and $j$-th factors and as the identity on all the other factors.

Consider the subbundle of conformal blocks
with fiber $\CB_{\zz}^\ell(N,\dots,N)\subset V^{\otimes |\lla|}$.
It is well known  \cite{KZ} that
this subbundle is invariant with respect to the KZ connection at level $\ell$.

For $\ell=1$ this subbundle is a line-bundle, with $P:=P_{\zz}(N,\ldots,N)$ being a nowhere zero section in it (see Theorem \ref{thm:P_basis}).
Since the KZ connection has regular singularities, the function
$$\tilde{P}:=P \cdot \prod_{1\leq i<j\leq mN} (z_i-z_j)^{-a_{ij}},$$
for some choice of numbers $a_{ij}$, must be a solution of the KZ equations with
$\ell=1$,
\begin{equation} \label{eq:KZ_form1}
\left( \frac{\partial}{\partial z_i}-
\frac{1}{m+1} \sum_{j\not= i} \frac{\Omega^{(i,j)}}{z_i-z_j}\right) \tilde{P}=0 ,
\qquad
i=1,\dots,|\lla|.
\end{equation}
Equivalent to (\ref{eq:KZ_form1}) is
\begin{equation} \label{eq:KZ_form2}
\left( \frac{\partial}{\partial z_i}- \frac{1}{m+1} \sum_{j\not= i} \frac{\Omega^{(i,j)}}{z_i-z_j}\right) P=P\cdot \sum_{j\not= i} \frac{a_{i,j}}{z_i-z_j}.
\end{equation}

In our conventions $\Omega^{(i,j)}$ reduces to the following operator:
$$\Omega^{(i,j)}(f) = f(y^{(k)}_i \leftrightarrow y^{(k)}_j\ \forall k).$$

\begin{theorem}
We have $a_{i,j}=-\frac{m}{m+1}$ for all $i\not= j$.
\end{theorem}

\begin{proof} Without loss of generality we choose $i=1$. The coefficient of $\prod_{k=1}^m Y^{(k)}_{ \{(k-1)N+1,\ldots,kN\}}$ on
the left hand side of (\ref{eq:KZ_form2})
is
\begin{equation}\label{eq:coeff}
-\sum_{j=N+1}^{mN} \frac{1}{R}\frac{1}{z_1-z_j} - \frac{1}{m+1}\left( \sum_{j=2}^N \frac{1}{z_1-z_j}\cdot\frac{1}{R} +
\sum_{j=N+1}^{mN} \frac{1}{z_1-z_j}\cdot \frac{1}{R^{1,j}} \right),
\end{equation}
where $R=R(z_1,\ldots,z_N| z_{N+1},\ldots,z_{2N}|\ldots|z_{(m-1)N+1},\ldots,z_{mN})$ is the denominator of the monomial $\prod_{k=1}^m Y^{(k)}_{ \{(k-1)N+1,\ldots,kN\}}$ in $P$, and $R^{1,j}=R(z_1\leftrightarrow z_j)$.
According to (\ref{eq:KZ_form2}) the expression (\ref{eq:coeff}) is equal to
$$\frac{1}{R} \sum_{j=2}^{mN} \frac{a_{1,j}}{z_1-z_j}.$$
Multiplying with $R$, and checking the residue at $z_{N+1}=z_1$ we obtain
$$-1-\frac{1}{m+1}(-1)=a_{1,N+1},$$
which implies $a_{1,N+1}=-m/(m+1)$. For the other $a_{i,j}$'s the result follows by symmetry.
\end{proof}

Summarizing, we obtained that
\begin{equation} \label{eq:tildeP}
\tilde{P}=P_{\zz}(N,\ldots,N)\cdot \prod_{1\leq i<j \leq mN} (z_i-z_j)^{\frac{m}{m+1}}
\end{equation}
satisfies the KZ differential equation (\ref{eq:KZ_form1}).

\begin{remark} In proving that (\ref{eq:tildeP}) satisfies the KZ equations we used a priori that
(\ref{eq:KZ_form2}) holds for some constants $a_{i,j}$. Equivalently,
  we made a particularly lucky choice of a hyperplane to take the
  residue of the coefficient of a well-chosen monomial. The fact that
  the resulting $\tilde{P}$ satisfies the KZ equations, equivalently,
  that other residues would give the same $a_{i,j}$'s encode
  remarkable identities. For example, one has
$$R \cdot \sum_{j=N+1}^{mN} \partial_{1,j}\left( \frac{1}{R} \right) =(1-m) \sum_{j=2}^N \frac{1}{z_1-z_j} + 2 \sum_{j=N+1}^{mN} \frac{1}{z_1-z_j}$$
(for the same $R$ as in the proof above). It would be interesting to see direct proofs of these properties; or even more interestingly, a combinatorial proof using the form (\ref{eq:concise}) of $P$.
\end{remark}

\begin{remark} We may consider the $\slm$ KZ equations
instead of the $\glm$ version (\ref{eq:KZ_form1}) above. The only
change in (\ref{eq:KZ_form1}) is that the $\glm$ Casimir operator has
to be replaced with the $\slm$ Casimir operator. In our conventions
the $\slm$ Casimir operator is
$$\Omega^{(i,j)}_{\slm}(f) = f(y^{(k)}_i \leftrightarrow
  y^{(k)}_j\ \forall k) -\frac{1}{m}f.$$ Arguments analogous to the
  above prove that
$$P_{\zz}(N,\ldots,N)\cdot \prod_{1\leq i<j \leq mN} (z_i-z_j)^{\frac{m-1}{m}}$$
is a solution to the $\slm$ KZ equations.
\end{remark}

\section{Products of level one conformal blocks are higher level conformal blocks}
\label{sec:products}

Recall that conformal blocks of level $1$ have dimension 1, namely
$$\dim \CB_{\zz}^1(\lla)=1, \qquad \text{for}\ d(\lla)=0\ \text{or}\ 1,$$
and that the functions $P_{\zz}(\lla)$ of Section \ref{section:def_P} form a basis in these spaces.
Recall also that for these $\lla$'s $e_{i,j}P_{\zz}(\lla)=0$ if $i<j$, as well as
\begin{equation} \label{building_blocks_properties}
\begin{aligned}[rl]
\left( e^{\zz}_{1,m} \right)^2 P_{\zz}(\lla)=0 & \qquad\ \text{if}\ d(\lla)=0, \\
       e^{\zz}_{1,m}           P_{\zz}(\lla)=0 & \qquad\ \text{if}\ d(\lla)=1.
\end{aligned}
\end{equation}

Now we are going to use these $P_{\zz}(\lla)$ functions as building
blocks to produce elements of $\CB_{\zz}^{\ell}(\lla)$ for
$\ell\geq 1$ and any $\lla$. In Section \ref{sect:span} we will show
that these elements span the space of conformal blocks for generic
$\zz$.

\begin{definition} Let $\lla\in \N^m$ be a partition, and let $k_1,\ldots,k_{|\lla|}$ be pairwise different integers.
We define $P^{(k_1,\ldots,k_{|\lla|})}(\lla)$ to be obtained from $P_{\zz}(\lla)$ by substituting
$$
\begin{aligned}[lcl]
y^{(j)}_{k_a} & \ \ \ \text{for}\ \ \  & y^{(j)}_a  & \ \ \ \text{and} \\
z_{k_a}       & \ \ \ \text{for}\ \ \  & z_a        &
\end{aligned}
$$
for all $j=1,\ldots,m$ and $a=1,\ldots,|\lla|$.
\end{definition}

For example
$$P^{(3,5)}(1,1)=\frac{y^{(1)}_3y^{(2)}_5}{z_3-z_5} +\frac{y^{(1)}_5y^{(2)}_3}{z_5-z_3}.$$

\begin{definition} \label{def:Q}
Let $\lla\in \N^m$ be a partition, and let $\ell\geq d(\lla)$. Put formally $\lambda_0=\ell+\lambda_m$. Let $\U=\{U_1,\ldots,U_{\ell}\}$ be a partitioning of the set $\{1,2,\ldots,|\lla|\}$ into $\ell$ subsets such that for every $j=0,\ldots,m-1$ we have
exactly $\lambda_j-\lambda_{j+1}$ of the $U_j$'s satisfying $|U_j|\equiv j \mod (m)$.
Consider each $U_j$ as an ordered set, ordered by the natural ordering of integers. Define the partition
$\mu^{(j)}\in \N^m$ to be the unique partition with $|\mu^{(j)}|=|U_j|$ and $d(\mu^{(j)})=0$ or $1$.
Define
$$Q(\U)=\prod_{j=1}^{\ell} P^{(U_j)}(\mu^{(j)}).$$
\end{definition}

\begin{proposition} \label{prop:Q_CB}
We have
$$Q(\U) \in \CB_{\zz}^{\ell}(\lla)$$
for any choice of $\U$.
\end{proposition}

\begin{proof} The proposition follows from the Leibnitz rule of differentiating products. Indeed, for $i<j$
$$e_{i,j} Q=\sum_{a=1}^{\ell} \left( e_{i,j} P^{(U_a)}(\mu^{(a)}) \cdot \prod_{b\not=a} P^{(U_b)}(\mu^{(b)})\right),$$
and the factor $e_{i,j} P^{(U_a)}(\mu^{(a)})$ in each term is 0.
Similarly,
$$\left(e^{\zz}_{1,m}\right)^{\ell+1-d(\lla)} = \sum \left(e^{\zz}_{1,m}\right)^{\bullet} P^{(U_1)}(\mu^{(1)}) \cdot  \left(e^{\zz}_{1,m}\right)^{\bullet}  P^{(U_2)}(\mu^{(2)})  \cdots   \left(e^{\zz}_{1,m}\right)^{\bullet}P^{(U_\ell)}(\mu^{(\ell)}),$$
where the $\bullet$'s stand for integers whose sum in each term is $\ell+1-d(\lla)$. We know that for $d(\lla)$ of the $\mu^{(b)}$ partitions
$e^{\zz}_{1,m} P^{(U_b)}(\mu^{(b)})=0$, and $\ell-d(\lla)$ of them $(e^{\zz}_{1,m})^2   P^{(U_b)}(\mu^{(b)})=0$, see (\ref{building_blocks_properties}). Hence, at least one of the factors in each term is 0.
\end{proof}

\begin{example}\rm
Let $m=2$, $\lla=(3,3)$, and $\ell=2$. Four of the 16 choices of $Q$ functions are
$$P^{(1,2,3,4,5,6)}(3,3), \qquad
P^{(1,2)}(1,1)\cdot P^{(3,4,5,6)}(2,2),$$
$$P^{(3,4)}(1,1)\cdot P^{(1,2,5,6)}(2,2),\qquad
P^{(5,6)}(1,1)\cdot P^{(1,2,3,4)}(2,2).$$
Calculation shows that these four polynomials form a basis for $\CB^2_{\zz}(3,3)$, for any choice of pairwise different $z_i$'s.
\end{example}

The dimension of $\CB_{\zz}^2(N,N)$ is $2^{N-1}$. Conjecturally the example above generalizes to the following.

\begin{conjecture} Let $m=2$. Let
$$\begin{aligned}
H= \{ \   \{U,V\}\ : & \ U,V \subset \{1,\ldots, 2N\}, U\cap V=\emptyset, U\cup V= \{1,\ldots, 2N\},\\
                 & \ (2i-1, 2i\in U \text{ or } 2i-1, 2i \in V) \ \text{for all } i=1,\ldots,N. \}
\end{aligned}$$
The collection of the $2^{N-1}$ polynomials (parameterized by $H$)
$$P^{(U)}\left(\frac{|U|}{2},\frac{|U|}2\right) \cdot P^{(V)}\left(\frac{|V|}{2},\frac{|V|}2\right)$$
 form a basis of the space of conformal blocks $\CB_{\zz}^2(N,N)$ (for any choice of pairwise different $z_i$'s).
\end{conjecture}


\section{Spanning the space of conformal blocks $\CB_{\zz}^{\ell}(\lla)$.}
\label{sect:span}

In this section we will define spaces $L^{\ell}(\lla)$ of certain
polynomials in the variables $y^{(j)}_a$ (without
$\zz$-dependence). Their relation with $Q(\U)$ functions, as well as a
recursive property of the $L^{\ell}(\lla)$ spaces will result in
proving that $\dim CB_{\zz}^{\ell}(\lla)=\dim L^{\ell}(\lla)$; as well
as in proving the fact that the $Q(\U)$ functions span the space of conformal
blocks for generic $\zz$. For notational simplicity in
Sections~\ref{sec:L_2}--\ref{sec:rec_L_2} we show the $m=2$ case in
detail, then sketch the necessary changes to obtain the results for
general $m$ in Section \ref{sec:conclusion}.

\subsection{The space $L^{\ell}(\lla)$ for $m=2$.} \label{sec:L_2}

For a subset $U=\{u_1<u_2<\ldots<u_r\}$ of $\N$ define the $R$-polynomial
$$R(U)=(y_{u_1}-y_{u_2})(y_{u_3}-y_{u_4})\cdots (y_{2\left[\frac{r}{2}\right]-1} - y_{2\left[\frac{r}{2}\right]}),$$
where $[x]$ is the integer part of $x$. For a collection of sets $U_1, \ldots, U_k$ the $R$-polynomial will be the product of the $R$-polynomials of the individual $U_i$'s. Below, the $U_i$ will always be disjoint sets. For example
$$R(\{1,2,3\},\{4,5\})=(y_1-y_2)(y_4-y_5).$$

\begin{definition} Let $\lla=(\lambda_1,\lambda_2)$ be a partition, and let $\ell\geq d(\lla)$. We define the vector space
$$
\begin{array}{lr}
L^{\ell}(\lla)=\spa\{ R(\U)\ |& \ \U \text{ is a partitioning of } \{1,\ldots,|\lla|\} \text{ into } \ell \text{ subsets}\ \ \\
& \text{(with some of the subsets possibly empty)},\ \  \\
& \text{such that } d(\lla) \text{ of them have odd cardinality},\ \  \\
& \ell-d(\lla) \text{ of them have even cardinality} \}.
\end{array}
$$
\end{definition}

\begin{example} \rm We have
$$\begin{array}{lll}
L^2(3,2)&=\spa\{  R(1,2,3,4,5), & \underbrace{R(\{1,2,3\},\{4,5\}), \ldots,R(\{4,5,6\},\{1,2\})}_{\binom{5}{3}},\\
&                               & \underbrace{R(\{1\},\{2,3,4,5\}), \ldots,R(\{5\},\{1,2,3,4\})}_{\binom{5}{1}} \} \\
&=\spa\{ (y_1-y_2)(y_3-y_4), & \underbrace{(y_1-y_2)(y_4-y_5),\ldots,(y_4-y_5)(y_1-y_2)}_{\binom{5}{3}}, \\
&                            & \underbrace{(y_2-y_3)(y_4-y_5), \ldots, (y_1-y_2)(y_3-y_4)}_{\binom{5}{1}} \},
\end{array}
$$ which turns out to be a 4-dimensional space. We have
$$\begin{array}{lll}
L^3(3,2)& = L^2(3,2)+ \spa\{R(\{1\},\{2,3\},\{4,5\}), \ldots\}\\
        & = L^2(3,2)+ \spa\{(y_2-y_3)(y_4-y_5),\ldots\}\\
        & = \spa\{ (y_i-y_j)(y_k-y_l): \#\{i,j,k,l\}=4\},
\end{array}
$$
which is 5-dimensional.
\end{example}

\subsection{Asymptotics of $Q(\U)$ functions for $m=2$.}

Let $m=2$ and $\lla=(\lambda_1, \lambda_2)$ be a partition. In this section we make the following notational simplifications: instead of the variables $y^{(1)}_a$ and $y^{(2)}_a$ we will write $x_a$ and $y_a$ respectively.

Define
$$
\begin{array}{lcl}
\alpha(\lla)=(0,1,1,2,2,3,3,\ldots,\lambda_2-1,\lambda_2-1,\underbrace{\lambda_2,\lambda_2,\ldots,\lambda_2}_{d(\lla)+1}) \in \N^{|\lla|} &
\text{if} & \lambda_2\not=0 \\
\alpha(\lla)=(0,\ldots,0) \in \N^{|\lla|} & \text{if} & \lambda_2=0.
\end{array}
$$

Straightforward calculation proves the following statement.

\begin{theorem} \label{thm:asym}
We have
$$\pm\lim_{z_{|\lla|}\to 0} \left( \lim_{z_{|\lla|-1}\to 0} \left( \ldots \lim_{z_1\to 0}\left( P_{\zz}(\lla) \cdot \prod_{i=1}^{|\lla|} z_i^{\alpha(\lla)_i}\right)\ldots \right) \right)=\ \ \ \ \ \ \ \ \ \ \ \ \ \ $$
$$\ \ \ \ \ \ \  (x_2y_1-x_1y_2)(x_4y_3-x_3y_4)\cdots(x_{2\lambda_2}y_{2\lambda_2-1}-x_{2\lambda_2-1}y_{2\lambda_2}) x_{2\lambda_2+1}\cdots x_{\lambda_1+\lambda_2}.$$
\end{theorem}

Recall the definition of $Q(\U)$ functions from Definition \ref{def:Q}. Theorem \ref{thm:asym} yields the following observation.

\begin{corollary}\label{cor:asym_Q_R}
For $1 \gg |z_{|\lla|}| \gg |z_{|\lla|-1}| \gg \ldots \gg |z_{1}| >0$ we have
$$Q(\U) \sim  \sum_{\substack{ \beta=(\beta_i)\\  \beta_i \geq \alpha_i}} A_{\beta}(\U) \prod_{i=1}^{|\lla|} z_i^{\beta_i}$$
where $\alpha_i$'s are suitable integers, and $A_{\beta}(\U)$'s are functions of $x_1,\ldots, x_{|\lla|}, y_1,\ldots,y_{|\lla|}$ such that
$A_\alpha(\U)|_{x_1=1, \ldots, x_{|\lla|}=1} = \pm R(\U)$.
\end{corollary}


\begin{corollary} \label{cor:dim_estimate}
For $\lla=(\lambda_1,\lambda_2)$, $\ell\geq d(\lla)$, and generic $\zz$ we have
$$\dim \CB_{\zz}^{\ell}(\lla) \geq \dim \spa \{ Q(\U): \U\} \geq \dim L^{\ell}(\lla).$$
\end{corollary}

\subsection{Recursion for the dimension of conformal blocks for $m=2$.}

\begin{theorem} \label{thm:recursion_CB}
Let $m=2$. The dimensions of $\CB_{\zz}^{\ell}(\lla)$ spaces satisfy the following recursions.
\begin{enumerate}
\item $\dim \CB^{\ell}(a,0)=1$, for $a=0,1,\ldots,\ell$.
\item $\dim \CB_{\zz}^{\ell}(a,b)=\dim \CB_{\zz}^{\ell}(a-1,b)+\dim \CB_{\zz}^{\ell}(a,b-1)$ for $0<a-b<\ell, b>0$.
\item $\dim \CB_{\zz}^{\ell}(a,a)=\dim \CB_{\zz}^{\ell}(a,a-1)$ for $a>0$.
\item $\dim \CB_{\zz}^{\ell}(b+\ell,b)=\dim \CB_{\zz}^{\ell}(b+\ell-1,b)$ for $b>0$.
\end{enumerate}
\end{theorem}

In other words, the dimensions satisfy the Pascal-triangle type of rule suggested in the following table for $\ell=3$.

\hspace{4 cm} {\bf Conformal blocks} \hspace{4 cm} {\bf dim}
\nobreak
$$
\begin{array}{cccccrccc}
 \CB_{\zz}^3(0,0) &        & & &\ \ \ \ \ \  & 1 & & & \\
&   \CB_{\zz}^3(1,0)          & & &\ \ \ \ \ \  & & 1 & &\\
 \CB_{\zz}^3(1,1) & &  \CB_{\zz}^3(2,0) & & \ \ \ \ \ \ & 1 & & 1 \\
&  \CB_{\zz}^3(2,1) & &  \CB_{\zz}^3(3,0) & &\ \ \ \ \ \  & 2 & & 1 \\
 \CB_{\zz}^3(2,2) & &  \CB_{\zz}^3(3,1) & &\ \ \ \ \ \  & 2 & & 3 \\
&  \CB_{\zz}^3(3,2) & &  \CB_{\zz}^3(4,1) & & \ \ \ \ \ \ & 5 & & 3 \\
 \CB_{\zz}^3(3,3) & &  \CB_{\zz}^3(4,2) & & \ \ \ \ \ \ & 5 & & 8 \\
&  \CB_{\zz}^3(4,3) & &  \CB_{\zz}^3(5,2) & \ \ \ \ \ \ & & 13 & & 8 \\
 \CB_{\zz}^3(4,4) & &  \CB_{\zz}^3(5,3) & & \ \ \ \ \ \ & 13 & & 21 \\
&  \CB_{\zz}^3(5,4) & &  \CB_{\zz}^3(6,3) & &\ \ \ \ \ \  & 34 & & 21 \\
\ldots & \ldots& \ldots &\ldots & \ \ \ \ \ \ & \ldots &\ldots & \ldots &\ldots\\
\end{array}
$$

\begin{proof} Dimensions of conformal blocks are
computed as structure constants of fusion rings, see
e.g. \cite{Z}, \cite{GV}. These structure constants satisfy natural recursion
relations, which in the case of ${\mathfrak{sl}_2}$ reduce exactly to the recursion
described in Theorem \ref{thm:recursion_CB}.
\end{proof}

\subsection{Recursion for the dimension of $L^{\ell}(\lla)$ spaces for $m=2$.}\label{sec:rec_L_2}

\begin{theorem} \label{thm:recursion_L}
Let $m=2$. The dimension of $L^{\ell}(\lla)$ spaces satisfy the
following properties.
\begin{enumerate}
\item $\dim L^{\ell}(a,0)=1$, for $a=0,1,\ldots,\ell$. \label{one}
\item $\dim L^{\ell}(a,b)\geq \dim L^{\ell}(a-1,b)+\dim L^{\ell}(a,b-1)$ for $0<a-b<\ell, b>0$. \label{two}
\item $\dim L^{\ell}(a,a)\geq \dim L^{\ell}(a,a-1)$ for $a>0$. \label{three}
\item $\dim L^{\ell}(b+\ell,b)\geq\dim L^{\ell}(b+\ell-1,b)$ for $b>0$. \label{four}
\end{enumerate}
\end{theorem}

\begin{proof} The space $L^{\ell}(a,0)$ is spanned by the constant 1 function. This proves (\ref{one}). To prove the other statements, we define linear maps
$$
\begin{array}{ll}
\Phi:L^{\ell}(\lla) \to L^{\ell}(\lambda_1+1,\lambda_2) & \text{ if } \ell>d(\lla) \\
\Psi:L^{\ell}(\lla) \to L^{\ell}(\lambda_1,\lambda_2+1) & \text{ if } d(\lla)>0,
\end{array}
$$
by
$$
\begin{array}{ll}
\Phi(g)=& g \ \ (\text{considered as a polynomial in one more variables}),\\
\Psi(g)=& SF\left(g\cdot(y_1+y_2+\ldots+y_{|\lla|} -d(\lla)y_{\lambda_1+\lambda_2+1})\right).
\end{array}
$$
Here $SF$ ($=$``square-free map'') is the linear map obtained by dropping the non-square-free terms.

\begin{example} For example, for the map $\Psi:L^3(3,1)\to L^3(3,2)$ we have
$$\begin{aligned}
R(\{1\},\{2,3\},\{4\})= (y_2-y_3) & \mapsto SF \left(  (y_2-y_3)(y_1+y_2+y_3+y_4-2y_5) \right)\\
                                  & = SF \left( (y_2-y_3)\left( (y_1-y_5)+ (y_4-y_5) \right) +(y_2-y_3)(y_2+y_3) \right)\\
                                  & = (y_2-y_3)(y_1-y_5)+(y_2-y_3)(y_4-y_5)\\
                                  & = R(\{2,3\},\{1,5\},\{4\})+ R(\{2,3\},\{4,5\},\{1\}).
\end{aligned}
$$
The map $\Phi$ is induced by the natural embedding
$$\C[y_1,\ldots,y_{\lambda_1+\lambda_2}] \to \C[y_1,\ldots,y_{\lambda_1+\lambda_2},y_{\lambda_1+\lambda_2+1}].$$
\end{example}

The map $\Phi$ is a linear embedding, provided we show that it indeed
maps into $L^{\ell}(\lambda_1+1,\lambda_2)$. It is enough to check
this for generators of $L^{\ell}(\lla)$. Let $R(\U)$ be such a
generator. Since $\ell>d(\lla)$ there is at least one even part in
$\U$. Adding the element $\lambda_1+\lambda_2+1$ to an even part we
define $\U'$. We have $\Phi(R(\U))=R(\U')$, showing that $\Phi$ indeed
maps to $L^{\ell}(\lambda_1+1,\lambda_2)$.

Now we need a lemma.

\begin{lemma} We have
\begin{equation} \label{eqn:simplified_Psi}
\Psi(R(\U))= R(\U) \cdot \left(\sum_{a} y_a - d(\lla)y_{\lambda_1+\lambda_2+1}\right),
\end{equation}
where the summation runs for those $a\in\{1,2,\ldots,|\lla|\}$ which are largest elements in odd parts of $\U$.
\end{lemma}

\begin{proof} The right hand side of (\ref{eqn:simplified_Psi}) is clearly square-free. The difference
$$R(\U)\cdot\left(y_1+y_2+\ldots+y_{|\lla|}
-d(\lla)y_{\lambda_1+\lambda_2+1}\right) - R(\U) \cdot \left(\sum_{a}
y_a - d(\lla)y_{\lambda_1+\lambda_2+1}\right)= R(\U)\cdot (\sum_b
y_b),$$ where the summation is for those $b$'s which are not largest
numbers of odd parts of $\U$. These $b$'s grouped in pairs occur in
$R(\U)$ in $(y_{b_1}-y_{b_2})$ factors. Using
$(y_{b_1}-y_{b_2})(y_{b_1}+y_{b_2})=y_{b_1}^2-y_{b_2}^2$ we obtain
that no term of the difference is square-free.  This proves the lemma.
\end{proof}

Our next claim is that $\Psi$ indeed maps into $L^{\ell}(\lambda_1,\lambda_2+1)$. Indeed, we have
$$\Psi(R(\U))=R(\U) \cdot \left(\sum_{a} y_a - d(\lla)y_{\lambda_1+\lambda_2+1}\right),$$
where $a$ is as in the lemma. This is further equal to
$$R(\U)\cdot \sum_a (y_a-y_{\lambda_1+\lambda_2+1})=\sum_{\U'} R(\U'),$$
where the sum runs for those $\U'$ that are obtained from $\U$ by adding $\lambda_1+\lambda_2+1$ to an odd part of $\U$.

Checking the coefficient of $y_{\lambda_1+\lambda_2+1}$ shows that $\Psi$ is an embedding.

Hence, we have that both $\Phi$ and $\Psi$ are linear embeddings. This proves (\ref{three}) and (\ref{four}). Our last claim is that the images of
$$\Phi:L^{\ell}(\lambda_1-1,\lambda_2) \to  L^{\ell}(\lambda_1,\lambda_2) \qquad\text{and}\qquad
\Psi:L^{\ell}(\lambda_1,\lambda_2-1) \to  L^{\ell}(\lambda_1,\lambda_2)$$
only intersect in 0 (for $d(\lla)<\ell$, $\lambda_2>0$). This follows from the fact that $\Psi(g)$ depends on the variable $y_{\lambda_1+\lambda_2+1}$ for all non-zero $g$, while no $\Phi(g)$ depends on this variable. This proves (\ref{two}).
\end{proof}

\begin{remark} \rm
The definition of the maps $\Phi$ and $\Psi$ are inspired by the definition of so-called iterated singular vectors from \cite{MV05}, \cite{RV}.
\end{remark}

\subsection{Conclusion.}\label{sec:conclusion}

Theorems \ref{thm:recursion_CB}, \ref{thm:recursion_L}, together with
Corollary \ref{cor:dim_estimate} prove inductively that for $m=2$ we
have
$$ \dim \CB_{\zz}^{\ell}(\lla)=\dim L^{\ell}(\lla).$$ Comparing with
the middle term of Corollary \ref{cor:dim_estimate} we obtain the
following theorem.

\begin{theorem}\label{thm:Q_generate}
Let $m=2$.
The $Q(\U)$ functions of Definition \ref{def:Q}, for all different choices of $\U$, span the space of conformal blocks $\CB_{\zz}^{\ell}(\lla)$ for generic $\zz$.
\end{theorem}

In fact Theorem \ref{thm:Q_generate} holds for $m\geq 2$ too. The proof goes along the same argument, as the one presented above for $m=2$, in Sections \ref{sec:L_2}--\ref{sec:rec_L_2}. The analogue of Theorem \ref{thm:asym} for higher $m$ is as follows.

\begin{proposition}
For $1 \gg |z_{|\lla|}| \gg |z_{|\lla|-1}| \gg \ldots \gg |z_{1}| >0$ we have
$$P_{\zz}(\lla) \sim  \sum_{\substack{ \beta=(\beta_i)\\  \beta_i \geq \alpha_i}} A_{\beta}(\lla) \prod_{i=1}^{|\lla|} z_i^{\beta_i}$$
where $\alpha_i$'s are suitable integers, and $A_{\beta}(\lla)$'s are functions of $y^{(i)}_a$ ($i=1,\ldots,m, a=1,\ldots,|\lla|$), such that
\begin{equation}\label{eqn:monster}
A_{\alpha}(\lla)=\pm\prod_{u=1}^m \prod_{v=1}^{\lambda_u-\lambda_{u+1}} \det\left(  y^{(i)}_{\sum_{k=u+2}^m \lambda_k + (u+1)\lambda_{u+1}+(v-1)u+j} \right)_{i,j=1,\ldots,u}.
\end{equation}
\end{proposition}


The $R(\U)$ functions for general $m$ are thus obtained by certain index shifts from (\ref{eqn:monster}), and the space $L^{\ell}(\lla)$ is obtained as the span of $R(\U)$ functions. After these definitions the proofs of Sections \ref{sec:L_2}--\ref{sec:rec_L_2} go through to the general $m\geq 2$ case.

\medskip

We conjecture that the $Q(\U)$ functions generate $\CB_{\zz}^{\ell}(\lla)$ not only for generic $\zz$, but for any $\zz=(z_1,\ldots,z_m)$ with $z_i\not= z_j$ for $i\not= j$.

\subsection{Bases in conformal blocks.}
A byproduct of our proof of Theorem \ref{thm:recursion_L}, is that we can choose a basis of $\CB_{\zz}^{\ell}(\lla)$ by just following how the relevant $\U$ partitions change at the possible paths of $\Phi$ and $\Psi$ maps. Here is what can be obtained this way, explained for $m=2$.

Fix $\lla=(\lambda_1,\lambda_2)$ and $\ell\geq \lambda_1-\lambda_2$. Consider a permutation $w_1,w_2,\ldots,w_{|\lla|}$ of the multiset $\{\underbrace{+,\ldots,+}_{\lambda_1},\underbrace{-,\ldots,-}_{\lambda_2}\}$ with the property
$$0\leq \sum_{i=1}^j w_i\leq \ell \qquad\text{for all}\ j.$$
Now we define a multiset $S_w^i$ of partitions of $\{1,\ldots,|\lla|\}$ into $\ell$ parts, for every $i=0,\ldots,|\lla|$, recursively:
\begin{itemize}
\item $S_w^0=\{\underbrace{\emptyset, \ldots, \emptyset}_{\ell}\}$;
\item if $w_{i+1}=+$, then  $S_w^{i+1}$ consists of partitions obtained from partitions in $S_w^i$ by adding the number $i+1$ to an even part;
\item if $w_{i+1}=-$, then  $S_w^{i+1}$ consists of partitions obtained from partitions in $S_w^i$ by adding the number $i+1$ to an odd part.
\end{itemize}
Define
$$Q_w=\sum_{\U\in S_w^{|\lla|}} Q(\U).$$
We obtained that the set
$$\{Q_w: w \text{ is a permutation of} \{\underbrace{+,\ldots,+}_{\lambda_1},\underbrace{-,\ldots,-}_{\lambda_2}\}, \text{with } 0\leq \sum_{i=1}^j w_i\leq \ell \text{ for all}\ j\}$$
is a basis of $\CB_{\zz}^{\ell}(\lla)$ for generic $\zz$.

\begin{example} \rm
For example, one of the 4 choices of $w$ for $\lla=(3,2)$ and $\ell=2$ is $++-+-$. Then we have
$$S_w^0=\{\emptyset, \emptyset\}, \ \ \ \  S_w^1=\{ \{ \{1\},\emptyset \} \},\ \ \ \  S_w^2=\{ \{ \{1\},\{2\} \} \},  $$
$$S_w^3=\{   \{ \{1,3\},\{2\} \}     ,    \{ \{1\},\{2,3\} \}\},  \ \ \ \  S_w^4=\{ \{ \{1,3,4\},\{2\} \} , \{ \{1\},\{2,3,4\} \}  \}, $$
$$S_w^5=\{   \{ \{1,3,4,5\},\{2\} \},
\{ \{1,3,4\},\{2,5\} \},
\{ \{1,5\},\{2,3,4\} \},
\{ \{1\},\{2,3,4,5\} \}   \}.
$$
Hence the corresponding basis vector is
$$Q_{++-+-}=Q(\{1,3,4,5\},\{2\})+ Q( \{1,3,4\},\{2,5\})+ Q(
\{1,5\},\{2,3,4\} )+ Q(\{1\},\{2,3,4,5\}).$$
Together with
$Q_{++--+}$, $Q_{+-++-}$, and $Q_{+-+-+}$
it forms a basis of $\CB^2_{\zz}(3,2)$ for generic $\zz$.
\end{example}

%

\section{Appendix: Decorated $P_{\zz}(\lla)$-functions}

In this section we present a generalization (a ``decorated version'') of the $P_{\zz}(\lla)$ function of Section~\ref{section:def_P}.

Let $m\geq 2$ be an integer and $\lla\in \N^m$ a partition as before.
Let the product of symmetric groups $S_{\lla}:=S_{\lambda_1}\times S_{\lambda_2}\times \ldots \times S_{\lambda_m}$ act on the polynomial ring
$\Z[z_1,\ldots,z_{|\lla|}]$ by permuting the variables
$$z_1,\ldots,z_{\lambda_1}\qquad\text{and}\qquad z_{\lambda_1+1},\ldots,z_{\lambda_1+\lambda_2},\qquad \text{and}\qquad \ldots \qquad \text{and}\qquad
z_{|\lla|-\lambda_m+1},\ldots,z_{|\lla|}$$
independently. A polynomial invariant under this action is called $\lla$-symmetric. A $\lla$-symmetric polynomial is unchanged if we permute e.g. the first $\lambda_1$ variables, or the next $\lambda_2$ variables, etc. Hence we may use the notation $h(z_{I_1},z_{I_2},\ldots,z_{I_m})$, if the $z_{I_j}$'s  are unordered sets of variables of cardinality $\lambda_j$.

\begin{definition}
Let $m\geq 2$, and let $\lla=(\lambda_1,\ldots,\lambda_m)\in \N^m$ be a partition. Let $h$ be a $\lla$-symmetric polynomial.
Recall that for a subset $U\subset \{1,\ldots,|\lla|\}$, we define $Y^{(j)}_U=\prod_{a\in U}y^{(j)}_a$. We define
$$P_{\zz}[h](\lla)=\sum_{\I} \frac{ h(z_{I_1},\ldots,z_{I_m})\cdot\prod_{j=1}^m Y_{I_j}^{(j)}}{R(z_{I_1}|z_{I_2}|\ldots|z_{I_m})},$$
where the summation runs for $\I=(I_1,\ldots,I_m)$ with $I_i\cap I_j=\emptyset$, $\cup I_j=\{1,\ldots,|\lla|\}$, $|I_j|=\lambda_j$.
\end{definition}

Clearly $P_{\zz}[1](\lla)=P_{\zz}(\lla)$ is the function defined in Section \ref{section:def_P}. Among the $P_{\zz}[h](\lla)$ functions one may often find several linearly independent conformal blocks.

\begin{definition} \label{def:schur} Let $Z_i^{(k)}$ be the $i$th elementary symmetric polynomial in the variables $z_1,\ldots,z_k$. For a partition $\mu=(\mu_1,\ldots,\mu_r)\in \N^r$ (for any $r\in \N$) let $$s^{(k)}(\mu)= \det \left( Z^{(k)}_{\mu_i+j-i} \right)_{r \times r}.$$
\end{definition}

For example
$$s^{(k)}(0)=1\ \text{for any}\ k,\qquad s^{(k)}(1)=z_1+\ldots+z_k,\qquad s^{(k)}(2)=z_1z_2+\ldots+z_{k-1}z_k,$$
$$s^{(k)}(1,1)=\det
\begin{pmatrix}
z_1+\ldots+z_k & z_1z_2+\ldots+z_{k-1}z_k\\
1 & z_1+\ldots+z_k
\end{pmatrix}.$$
The polynomial $s^{(k)}(\mu)$ is a symmetric polynomial in $z_1,\ldots,z_k$, and does not depend on any $z_a$ variable for $a>k$. Hence it is a $\lla$-symmetric function for all $\lla$ with $\lambda_1=k$.

\begin{example} \rm
Calculation shows that
$$P_{\zz}[h](4,2) \in \CB_{\zz}^2(4,2) \qquad\text{for}\qquad h=s^{(4)}(0), s^{(4)}(1), s^{(4)}(1,1).$$
Moreover, these polynomials are linearly independent for generic $\zz$. We also have
$$P_{\zz}[h](4,2) \in \CB_{\zz}^3(4,2) \qquad\text{for}\qquad h=s^{(4)}(0), s^{(4)}(1), s^{(4)}(1,1),$$
$$\phantom{P_{\zz}[h](\lla) \in \CB_{\zz}^3(4,2) \qquad\text{for}\qquad h=   }  \ \ \      s^{(4)}(2), s^{(4)}(2,1), s^{(4)}(2,2),$$
and these polynomials are linearly independent for generic $\zz$.  It is known that
$\dim \CB_{\zz}^2(4,2)$ $=$ $4$, $\dim \CB_{\zz}^3(4,2)=8$. Hence our $P_{\zz}[h](4,2)$ functions above do not span the relevant spaces of conformal blocks. Moreover, for the 9-dimensional spaces $\CB^{>3}(4,2)$ no new linearly independent $P_{\zz}[h](4,2)$ function can be found.
\end{example}

The subspace of conformal blocks spanned by $P_{\zz}[h]$ functions, and the relation to Schubert calculus is subject to future study.


\begin{thebibliography}{Mat90}

\bibitem[AB]{ab}
M. F. Atiyah, R. Bott.
\newblock The moment map and equivariant cohomology.
\newblock {\em Topology, Vol. 23, Issue 3, 1--28, 1984}

\bibitem[FSV1]{FSV1}
B. Feigin, V. Schechtman, A. Varchenko.
\newblock On algebraic equations satisfied by hypergeometric correlators in WZW models. I,
\newblock {\em Comm. Math. Phys. 163 (1994), 173--184}

\bibitem[FSV2]{FSV2}
B. Feigin, V. Schechtman, A. Varchenko.
\newblock On algebraic equations satisfied by hypergeometric correlators in WZW models. II,
\newblock {\em Comm. in Math. Phys. v. 170, No. 1, (1994) 219--247; math.hep-th/9407010}

\bibitem[FSV]{FSV}
G. Felder, L. Stevens, A. Varchenko.
\newblock Elliptic Selberg integrals and conformal blocks.
\newblock {\em Math. Res. Lett.  10 (2003), no. 5-6, 671--684}

\bibitem[GV]{GV}
S. Gusein-Zade, A. Varchenko.
\newblock Verlinde algebras and the intersection form on vanishing cycles.
\newblock {\em Selecta Math. (N.S.) 3 (1997), no. 1, 79--97}

\bibitem[KL]{KL}
D. Kazhdan, G. Lusztig.
\newblock Tensor categories arising from affine Lie algebras I-V.
\newblock {\em J. Amer. Math. Soc., 6(1993), 905--947; ibid., 949--1011; 7(1994), 335--381; ibid., 383--454}

\bibitem[KZ]{KZ}
V.G. Knizhnik, A.B. Zamolodchikov.
\newblock Current Algebra and Wess-Zumino Model in Two-Dimensions.
\newblock {\em Nucl. Phys. B 247 (1984), 83--103}

\bibitem[L]{lascoux}
A. Lascoux.
\newblock Symmetric Functions and Combinatorial Operations on Polynomials.
\newblock {\em CMBS 99, AMS, 2003}.

\bibitem[LV]{LV}
E. Looijenga, A. Varchenko.
\newblock Unitarity of ${\rm SL} (2)$-conformal blocks in genus zero.
\newblock {\em arXiv:0810.4310, [v1], 1--15}

\bibitem[MV]{MV}
E. Mukhin, A. Varchenko.
\newblock Remarks on critical points of phase functions and norms of Bethe vectors;
\newblock {\em Arrangements---Tokyo 1998, 239--246, Adv. Stud. Pure Math., 27, Kinokuniya, Tokyo, 2000}

\bibitem[MV05]{MV05}
E. Mukhin, A. Varchenko.
\newblock Norm of a Bethe vector and the Hessian of the master function.
\newblock {\em Compos. Math. 141 (2005), no. 4, 1012--1028.}


\bibitem[R]{R}
Ramadas, T. R. 
\newblock The ``Harder-Narasimhan trace'' and unitarity of the KZ/Hitchin connection: genus 0. 
\newblock {\em Ann. of Math. (2) 169 (2009), no. 1, 1--39.}

\bibitem[RV]{RV}
N. Reshetikhin, A. Varchenko.
\newblock Quasiclassical asymptotics of solutions to the KZ equations.
\newblock {\em  Geometry, topology, \& physics, 293--322, Conf. Proc. Lecture Notes Geom. Topology, IV, Int. Press, Cambridge, MA, 1995.}




\bibitem[V08]{anv_related} A. Varchenko.  \newblock A Selberg Integral
  Type Formula for an $sl_2$ One-Dimensional Space of Conformal
  Blocks.  \newblock {\em arXiv:0810.3355, 2008}

\bibitem[Z]{Z}
J.-B.\,Zuber.
\newblock Graphs and reflection groups.
\newblock {\em Comm. Math. Phys. 179 (1996), no. 2, 265--294}


\end{thebibliography}
\end{document}